\documentclass[8pt]{article}
\usepackage{amsmath}
\usepackage{amssymb}
\usepackage{amsfonts}
\usepackage{amssymb}
\usepackage{graphicx}
\usepackage{lscape}
\renewcommand{\baselinestretch}{1.75}
\newtheorem{theo}{Theorem}
\newtheorem{lem}{Lemma}

\def\twoem{\,\,\,\,\,\,}

\def\var{\mathrm{var}}

\newcounter{hyp}

\begin{document}

\begin{titlepage}

\title{Propagation of Memory Parameter from Durations to Counts}

\renewcommand{\baselinestretch}{1.2}

\author{ Rohit Deo\thanks{New York University, 44 W. 4'th Street, New York NY 10012 USA} \and Clifford M.
  Hurvich$^*$ \and Philippe Soulier\thanks{Universit\'e Paris X, 200 avenue de la R\'epublique, 92001 Nanterre cedex, France
  \newline The authors thank Xiaohong Chen and Raymond Thomas for helpful comments and suggestions.} \and Yi
  Wang$^*$}
%\date{}

\maketitle

\begin{abstract}
We establish sufficient conditions on durations that are
stationary with finite variance and memory parameter $d \in
[0,1/2)$ to ensure that the corresponding counting process $N(t)$
satisfies $\textmd{Var} \, N(t) \sim C t^{2d+1}$ ($C>0$) as $t
\rightarrow \infty$, with the same memory parameter $d \in
[0,1/2)$ that was assumed for the durations. Thus, these
conditions ensure that the memory in durations propagates to the
same memory parameter in counts and therefore in realized
volatility. We then show that any Autoregressive Conditional
Duration ACD(1,1) model with a sufficient number of finite moments
yields short memory in counts, while any Long Memory Stochastic
Duration model with $d>0$ and all finite moments yields long
memory in counts, with the same $d$. Finally, we present a result
implying that the only way for a series of counts aggregated over
a long time period to have nontrivial autocorrelation is for the
short-term counts to have long memory. In other words, aggregation
ultimately destroys all autocorrelation in counts, if and only if
the counts have short memory. \newline \newline \textit{KEYWORDS}:
Long Memory Stochastic Duration, Autoregressive Conditional
Duration, Rosenthal-type Inequality.

\end{abstract}

\end{titlepage}

\newpage

\section{Introduction}
There is a growing literature on long memory in volatility of
financial time series. See, e.g., Robinson (1991), Bollerslev and
Mikkelsen (1996), Robinson and Henry (1999), Deo and Hurvich
(2001), Hurvich, Moulines and Soulier (2005). Long memory in
volatility, which has been repeatedly found in the empirical
literature, plays a key role in the forecasting of realized
volatility (Andersen, Bollerslev, Diebold and Labys 2001, Deo,
Hurvich and Lu 2005), and has important implications on option
pricing (see Comte and Renault 1998).

Given the increasing availability of transaction-level data it is
of interest to explain phenomena observed at longer time scales
from equally-spaced returns in terms of more fundamental
properties at the transaction level. Engle and Russell (1998)
proposed the Autoregressive Conditional Duration (ACD) model to
describe the durations between trades, and briefly explored the
implications of this model on volatility of returns in discrete
time, though they did not determine the persistence of this
volatility, as measured, say, by the decay rate of the
autocorrelations of the squared returns. Deo, Hsieh and Hurvich
(2005) proposed the Long-Memory Stochastic Duration (LMSD) model,
and began an empirical and theoretical exploration of the question
as to which properties of durations lead to long memory in
volatility, though the theoretical results presented there were
not definitive.

The collection of time points $\cdots t_{-1} < t_0 \leq 0 < t_1 <
t_2 < \cdots$ at which a transaction (say, a trade of a particular
stock on a specific market) takes place, comprises a point process,
a fact which was exploited by Engle and Russell (1988). These event
times $\{t_k\}$ determine a \textit{counting process},
\[
N(t)= Number \,\, of \,\, Events \,\, in \,\, (0,t].
\]
For any fixed time spacing $\Delta t >0$, one can define the
\textit{counts} $\Delta N_{t^\prime} = N(t^\prime \Delta t ) -
N((t^\prime - 1) \Delta t)$, the number of events in the
$t^\prime$'th time interval of width $\Delta t$, where $t^\prime =
1,2,\cdots$. The event times $\{t_k\}_{k=-\infty}^{\infty}$ also
determine the \textit{durations}, given by
$\{\tau_k\}_{k=-\infty}^{\infty}$, $\tau_k = t_k - t_{k-1}$.

Both the ACD and LMSD models imply that the doubly infinite sequence
of durations $\{\tau_k\}_{k=-\infty}^{\infty}$ are a stationary time
series, i.e., there exists a probability measure $P^0$ under which
the joint distribution of any subcollection of the $\{\tau_k\}$
depends only on the lags between the entries. On the other hand, a
point process $N$ on the real line is stationary under the measure
$P$ if $P(N(A))=P(N(A+c))$ for all real $c$. A fundamental fact
about point processes is that in general (a notable exception is the
Poisson process) there is no single measure under which both the
point process $N$ and the durations $\{\tau_k\}$ are stationary,
i.e., in general $P$ and $P^0$ are not the same. Nevertheless, there
is a one-to-one correspondence between the class of measures $P^0$
that determine a stationary duration sequence and the class of
measures $P$ that determine a stationary point process. The measure
$P^0$ corresponding to $P$ is called the \textit{Palm distribution}.
The counts are stationary under $P$, while the durations are
stationary under $P^0$.

Deo, Hsieh and Hurvich (2005) pointed out, using a theorem of
Daley, Rolski and Vesilo (2000) that if durations are generated by
an $ACD$ model and if the durations have tail index $\kappa \in
(1,2)$ under $P^0$, then the resulting counting process $N(t)$ has
long range count dependence with memory parameter $d \geq
1-\kappa/2$, in the sense that $\textmd{Var} \, N(t) \sim C
n^{1+2d}$ ($C>0$) as $t \rightarrow \infty$, under $P$. This,
together with the model for returns at equally spaced intervals of
time given in Deo, Hsieh and Hurvich (2005) implies that realized
volatility has long memory in the sense that the $n$-term partial
sum of realized volatility has a variance that scales as $C_2
n^{2d+1}$ as $n\rightarrow \infty$, where $C_2 >0$. Deo, Hsieh and
Hurvich (2005) also showed that if durations are generated by an
LMSD model with memory parameter $d$ under $P^0$ then counts have
long memory with memory parameter $d^{counts} \geq d$, but
unfortunately this conclusion was established only under the
duration-stationary measure $P^0$, and not under the
count-stationary measure $P$. This gap can be bridged using
methods described in this paper. Still, the results we have
described above merely give lower bounds for the memory parameter
in counts.

In this paper, we will establish sufficient conditions on durations
that are stationary with finite variance and memory parameter $d \in
[0,1/2)$ under $P^0$ to ensure that the corresponding counting
process $N(t)$ satisfies $\textmd{Var} \, N(t) \sim C t^{2d+1}$
($C>0$) as $t \rightarrow \infty$ under $P$, with the same memory
parameter $d \in [0,1/2)$ that was assumed for the durations. Thus,
these conditions ensure that the memory in durations propagates to
the same memory parameter in counts and therefore in realized
volatility.

Next, we will verify that the sufficient conditions of our Theorem
\ref{theo:general} are satisfied for the ACD(1,1) model assuming
finite $8+\delta$ moment ($\delta>0$) of the durations under
$P^0$, and for the LMSD model with any $d \in [0,1/2)$ assuming
that the multiplying shocks have all moments finite. Thus, any
ACD(1,1) model with a sufficient number of finite moments yields
short memory in counts, while any LMSD model with $d>0$ and all
finite moments yields long memory in counts. These results for the
LMSD and ACD(1,1) models are given in Theorems \ref{theo:LMSD} and
\ref{theo:ACD}, respectively. Lemma \ref{lem:partialsumLMSD},
which is used in proving Theorem \ref{theo:LMSD}, provides a
Rosenthal-type inequality for moments of absolute standardized
partial sums of durations under the LMSD model, and is of interest
in its own right.

Finally, we present a result (Theorem \ref{theo:Autocor}) implying
that if counts have memory parameter $d \in [0, 1/2)$ then further
aggregations of these counts to longer time intervals will have a
lag-1 autocorrelation that tends to $2^{2d}-1$ as the level of
aggregation grows. Interestingly, this limit is zero if and only
if $d=0$. Thus, one of the important functions of long memory in
counts is that it allows the counts to have a non-vanishing
autocorrelation even as $\Delta t$ grows, as was found by Deo,
Hsieh and Hurvich (2005) to occur in empirical data. By contrast,
short memory in counts implies that counts at long time scales
(large $\Delta t$) are essentially uncorrelated, in contradiction
to what is seen in actual data. To summarize, aggregation
ultimately destroys all autocorrelation in counts, if and only if
the counts have short memory.

\section{Theorems on the propagation of the memory parameter}

Let $E$, $E^0$, $\textmd{Var}$, $\textmd{Var}^0$ denote
expectations and variances under $P$ and $P^0$, respectively.
Define $\mu=E^0(\tau_k)$ and $\lambda=\frac{1}{\mu}.$ Our main
theorem uses the assumption that $P^0$ is $\{\tau_k\}$-mixing,
defined as follows. Let $\mathcal N =
\sigma(\{\tau_k\}_{k=-\infty}^{\infty})$ and $\mathcal F_n=
\sigma(\{\tau_k\}_{k=n}^{\infty})$. We say that $P^0$ is
$\{\tau_k\}$-mixing if
\[
\lim_{n\rightarrow \infty} \sup_{B \in \mathcal N \cap \mathcal
F_n} |P^0(A \cap B) - P^0(A) P^0(B)| = 0 \] for all $A \in
\mathcal N$.

\begin{theo}\label{theo:general}

Let $\{ \tau_k \}$ be a duration process such that the following
conditions hold:

$i)$ $\{ \tau_k \}$ is stationary under $P^0$.

$ii)$ $P^0$ is $\{ \tau_k \}$-mixing.

$iii)$ $\exists ~ d \in [0, \frac{1}{2})$ such that
\[
Y_n (s) = \frac{\sum_{k=1}^{\lfloor ns \rfloor} ( \tau_k - \mu)
}{n^{1/2+d}}, ~~~~~ s \in [0,1]
\]
converges weakly to $\sigma B_{1/2+d}(\cdot)$ under $P^0$, where
$\sigma>0$ and $B_{1/2+d} (\cdot)$ is fractional Brownian motion
if $0 < d < \frac{1}{2}$ or standard Brownian motion $B_{1/2}=B$
if $d=0$.

$iv)$
\[
\sup_n E^0 \Big|\frac{\sum_{k=1}^n ( \tau_k - \mu)
}{n^{1/2+d}}\Big|^p < \infty ~~~~~ \left\{ \begin{array} {l}
for ~ all ~ p>0, ~ if ~ d \in (0,\frac{1}{2}) \\
for ~ p=8+\delta, ~ \delta > 0, ~ if ~ d=0 ~~~~~ \cdot
\end{array} \right.
\]

Then the induced counting process $N(t)$ satisfies $\textmd{Var}N(t)
\sim Ct^{2d+1}$ under $P$ as $t \rightarrow \infty$ where $C>0$.
\end{theo}
\textbf{Remark:} Inspection of the proof of Theorem
\ref{theo:general} reveals that if $d>0$, only $4/(0.5-d)+\delta$
finite moments are needed, where $\delta>0$ is arbitrarily small.
The closer $d$ is to $1/2$, the larger the number of finite
moments required.
\newline
\newline
\textbf{Remark:} As pointed out by Nieuwenhuis (1989), if
$\{\tau_k\}$ is strong mixing under $P^0$ then $P^0$ is
$\{\tau_k\}$-mixing. This weaker form of mixing is essential for
our purposes since even Gaussian long-memory processes are not
strong mixing. See Gu\'egan and Ladoucette (2001).

\subsection{LMSD Process}

Define the LMSD process $\{\tau_k\}_{k=-\infty}^{\infty}$ for $d
\in [0,\frac{1}{2})$ as
\[
\tau_k = e^{h_k} \epsilon_k
\]
where under $P^0$ the $\epsilon_k \ge 0$ are $i.i.d.$ with all
moments finite, and $h_k=\sum_{j=0}^\infty b_j e_{k-j}$, the $\{
e_k \}$ are \textit{i.i.d.} Gaussian with zero mean, independent
of $\{ \epsilon_k \}$, and
\[
b_j \sim \left \{ \begin{array} {l} Cj^{d-1} ~~~~~~~~~~~~ if~d \in (0,\frac{1}{2}) \\
Ca^j,\,\,|a|<1 ~~~ ~if~d=0
\end{array} \right.
\]
($C\neq 0$) as $j \rightarrow \infty$. Note that for convenience,
we nest the short-memory case ($d=0$) within the LMSD model, so
that the allowable values for $d$ in this model are $0 \leq d <
1/2$.

\begin{theo}\label{theo:LMSD}
If the durations $\{ \tau_k \}$ are generated by the LMSD process
with $d\in[0,1/2)$, then the induced counting process $N(t)$
satisfies $\textmd{Var}N(t) \sim Ct^{2d+1}$ under $P$ as $t
\rightarrow \infty$ where $C>0$.
\end{theo}
To establish Theorem \ref{theo:LMSD}, we will use the following
Rosenthal-type inequality.

\begin{lem}\label{lem:partialsumLMSD}
For durations $\{ \tau_k \}$ generated by the LMSD process with $d
\in [0,\frac{1}{2})$, for any fixed positive integer $p \ge 2$,
$E^0\{ |y_n-E^0(y_n)|^p \}$ is bounded uniformly in $n$, where
\[
y_n = \frac{\sum_{k=1}^n \tau_k}{n^{1/2+d}}~~~~~.
\]
\end{lem}

\subsection{ACD(1,1) Process}

Define the ACD(1,1) process $\{\tau_k\}_{k=-\infty}^{\infty}$ as
\begin{eqnarray*}
\tau_k &=& \psi_k \epsilon_k \\
\psi_k &=& \omega + \alpha \tau_{k-1} + \beta \psi_{k-1}
\end{eqnarray*}
with $\omega>0,\alpha > 0, \beta \ge 0$ and $\alpha+\beta < 1$,
where under $P^0$, $\epsilon_k \ge 0$ are \textit{i.i.d.} with
mean 1. We will assume further that under $P^0$, $\epsilon_k$ has
a density $g_\epsilon$ such that $\int_0^\theta g_\epsilon(x) dx
>0, \forall ~ \theta >0$
and $E^0(\tau_k^{8+\delta}) < \infty$ for some $\delta>0$.

Nelson (1990) guarantees the existence of the doubly-infinite
ACD(1,1) process $\{\tau_k\}_{k=-\infty}^{\infty}$, which in our
terminology is stationary under $P^0$.

\begin{theo}\label{theo:ACD}
Suppose that the durations $\{ \tau_k \}$ are generated by the
ACD(1,1) model, with the additional assumptions stated above.
Then the induced counting process $N(t)$ satisfies
$\textmd{Var}N(t) \sim Ct$ under $P$ as $t \rightarrow \infty$
where $C>0$.
\end{theo}

\section{Autocorrelation of Aggregated Counts}

\begin{theo} \label{theo:Autocor}
Let $\{X_t\}$ be a stationary process such that
$\textmd{Var}(\sum_{t=1}^n X_t) \sim Cn^{1+2d}$ as $n \rightarrow
\infty$, where $C \neq 0$ and $d \in [0,1/2)$. Then
\[
\lim_{n\rightarrow \infty} Corr \left [ \sum_{t=1}^n X_t ,
\sum_{t=n+1}^{2n} X_t \right ] = 2^{2d}-1.
\]
\end{theo}

\textbf{Proof:}
\[
\textmd{Var} \left [ \sum_{t=1}^{2n} X_t \right ] = 2\,
\textmd{Var} \left [ \sum_{t=1}^{n} X_t \right ] + 2
\,\textmd{Cov} \left [ \sum_{t=1}^n X_t , \sum_{t=n+1}^{2n} X_t
\right ].
\]

Thus,
\[
\textmd{Cov} \left [ \sum_{t=1}^n X_t , \sum_{t=n+1}^{2n} X_t \right
]= .5 \left ( \textmd{Var} \left [ \sum_{t=1}^{2n} X_t \right ] - 2
\textmd{Var} \left [ \sum_{t=1}^{n} X_t \right ] \right ).
\]
The result follows by noting that $\lim_{n \rightarrow \infty}
n^{-2d-1} \textmd{Var} (\sum_{t=1}^n X_t) = C$, where $C \neq 0$.
$\Box$

This theorem has an interesting practical interpretation. If we
write $X_k = N[k\Delta t] - N[(k-1) \Delta t ]$ where $\Delta t>0$
is fixed, then $X_k$ represents the number of events (count) in a
time interval of width $\Delta t$, e.g. one minute. Thus,
$\sum_{k=1}^n X_k$ is the number of events in a time interval of
length $n$ minutes, e.g. one day. The theorem implies that as the
level of aggregation ($n$) increases, the lag-1 autocorrelation of
the aggregated counts will approach a nonzero constant if and only
if the non-aggregated count series $\{X_k\}$ has long memory. In
other words, the only way for a series of counts over a long time
period to have nontrivial autocorrelation is for the short-term
counts to have long memory. Since in practice long-term counts do
have substantial autocorrelation (see Deo, Hsieh and Hurvich
2005), it is important to use only the models for durations that
imply long memory in the counting process (LRcD). Examples of such
models include the LMSD model (see Theorem \ref{theo:LMSD}), and
ACD models with infinite variance (see Daley, Rolski and Vesilo
(2000), and Theorem 2 of Deo, Hsieh and Hurvich, 2005).

\section{Appendix: Proofs}

Let $P$ denote the stationary distribution of the point process $N$
on the real line, and let $P^0$ denote the corresponding Palm
distribution. $P$ determines and is completely determined by the
stationary distribution $P^0$ of the doubly infinite sequence
$\{\tau_k \}_{k=-\infty}^\infty$ of durations. Note that the
counting process $N$ is stationary under $P$, the durations are
stationary under $P_0$, but in general there is no single
distribution under which both the counting process and the durations
are stationary. For more details on the correspondence between $P$
and $P^0$, see Daley and Vere-Jones (2003), Baccelli and
Br$\acute{\text{e}}$maud (2003), or Nieuwenhuis (1989).

Following the standard notation for point processes on the real line
(see, e.g., Nieuwenhuis 1989, p. 594), we assume that the event
times $\{t_k\}_{k=-\infty}^{\infty}$ satisfy
\[
\ldots < t_{-1} < t_0 \leq 0 < t_1 < t_2 < \ldots .
\]

Let
\[
u_k = \left \{ \begin{array} {l} t_1 ~~~~~ if ~ k=1 \\
\tau_k ~~~~~ if ~ k \ge 2 ~~~~~ \cdot
\end{array} \right.
\]
Here, the random variable $t_1>0$ is the time of occurrence of the
first event following $t=0$. For $t > 0$, define the count on the
interval $(0,t]$, $N(t):=N(0,t]$, by
\begin{eqnarray*}
N(t) &=& \max\{ s:\sum_{i=1}^s u_i \le t \},~~~u_1 \le t \\ &=& 0,
~~~~~~~~~~~~~~~~~~~~~~~~~~~u_1>t.
\end{eqnarray*}

Throughout the paper, the symbol $\Longrightarrow$ denotes weak
convergence in the space $D[0,1]$.

\textbf{Proof of Theorem \ref{theo:general}:}

By assumption $iii)$, $Y_n \Longrightarrow \sigma B_{1/2+d}$ under
$P^0$, where $\sigma>0$. First, we will apply Theorem 6.3 of
Nieuwenhuis (1989) to the durations
$\{\tau_k\}_{k=-\infty}^{\infty}$ to conclude that $Y_n
\Longrightarrow \sigma B_{1/2+d}$ under $P$. Since the
$\{\tau_k\}_{k=-\infty}^{\infty}$ are stationary under $P^0$ and
are generated by the shift to the first event following time zero
(see Nieuwenhuis 1989, p. 600), and since we have assumed that
$P^0$ is $\{\tau_k\}$-mixing, his Theorem 6.3 applies. It follows
that $Y_n \Longrightarrow \sigma B_{1/2+d}$ under $P$. We next
show that the suitably normalized counting process converges to
the same limit under $P$.

Define
\[
\tilde Y_n(s)=\frac {\sum_{k=1} ^{\lfloor ns
\rfloor}(u_k-\mu)}{n^{1/2+d}} \twoem, \twoem s \in [0,1].
\]
Note that for all $s$, $\tilde Y_n(s) = Y_n(s) +
n^{-(1/2+d)}(u_1-\tau_1)$. From Baccelli and
Br$\acute{\text{e}}$maud (2003, Equation 1.4.2, page 33), for any
measurable function $h$,
\begin{equation} \label{eq:ETau1}
E[h(\tau_1)] = \lambda E^0 [\tau_1 h(\tau_1)] \twoem.
\end{equation}
Since $u_1 \leq \tau_1$, and since assumption $iv)$ implies that
$\tau_1$ has finite variance under $P^0$, using $h(x)=x$ in
(\ref{eq:ETau1}), it follows that $n^{-(1/2+d)}(u_1-\tau_1)$ is
$o_p(1)$ under $P$. Thus, $\tilde Y_n \Longrightarrow \sigma
B_{1/2+d}$ under $P$.

Let
\begin{equation}\label{eq:defZt}
Z(t)=\frac{N(t)-t/\mu}{t^{1/2+d}} ~~~.
\end{equation}

By Iglehart and Whitt (1971, Theorem 1), it follows that $Z(t)
\stackrel d \rightarrow \tilde C B_{1/2+d}(1)$ under $P$ as $t
\rightarrow \infty$, where $\tilde C>0$. Furthermore, by Lemma
\ref{lem:supz}, $Z^2(t)$ is uniformly integrable under $P$ and hence
$\lim_t \textmd{Var} [Z(t)]=\tilde C^2 \textmd{Var} [B_{1/2+d}(1)]$.
The theorem is proved. $\Box$

\textbf{Proof of Theorem \ref{theo:LMSD}:}

We simply verify that the conditions of Theorem \ref{theo:general}
hold for this process.

By definition $\{ \tau_k \}$ is stationary under $P^0$ and by
Lemma \ref{lem:mixing}, $P^0$ is $\{ \tau_k \}$ mixing. By
Surgailis and Viano (2002), $Y_n \Longrightarrow \sigma B_{1/2+d}$
under $P^0$, where $\sigma>0$ and by Lemma
\ref{lem:partialsumLMSD}, $ \sup_n E^0 \Big|\frac{\sum_{k=1}^n (
\tau_k - \mu) }{n^{1/2+d}}\Big|^p < \infty$ for all $p$. Thus, the
result is proved. $\Box$

\textbf{Proof of Theorem \ref{theo:ACD}:}

We simply verify that the conditions of Theorem \ref{theo:general}
hold for this process.

By Lemma \ref{lem:mixing}, $\{ \tau_k \}$ is exponential
$\alpha$-mixing, and hence strong mixing and thus by Nieuwenhuis
(1989), $P^0$ is $\{\tau_k\}$-mixing. Furthermore, since all
moments of $\tau_k$ exist up to order $8+\delta,\delta>0$, we can
apply results from Doukhan (1994) to obtain
\begin{equation}\label{eq:convergence}
Y_n \Rightarrow C B,
\end{equation}
if $\frac{1}{n} \var(\sum_{k=1}^n \tau_k) \rightarrow C^2
>0$, as $n \rightarrow \infty$.

It is well known that the GARCH(1,1) model can be represented as an
ARMA(1,1) model, see Tsay (2002). Similarly, the ACD(1,1) model can
also be re-formulated as an ARMA(1,1) model,
\begin{equation}\label{eq:arma}
\tau_k = \omega + (\alpha+\beta) \tau_{k-1} + (\eta_k - \beta
\eta_{k-1})
\end{equation}
where $\eta_k = \tau_k - \psi_k$ is white noise with finite variance
since $E(\tau_k^{8+\delta}) < \infty$. The autoregressive and moving
average parameters of the resulting ARMA(1,1) model are
$(\alpha+\beta)$ and $\beta$, respectively.

It is also known that for any stationary invertible ARMA model $\{
z_k \}$, $n \var(\bar z) \rightarrow 2 \pi f_z(0)$, where $f_z(0)$
is the spectral density of $\{ z_k \}$ at zero frequency. For an
ARMA(1,1) process, $f_z(0)>0$ if the moving average coefficient is
less than 1. Here, since $0 \le \beta <1$, we obtain $\frac{1}{n}
\var(\sum_{k=1}^n \tau_k) = n \var(\bar \tau) \rightarrow 2 \pi
f_\tau(0) > 0$, as $n \rightarrow \infty$. Therefore
(\ref{eq:convergence}) follows.

Define $ y_n = \frac{1}{\sqrt{n}} \sum_{k=1}^n \tau_k$. Since all
moments of $\tau_k$ are bounded up to order $8+\delta$,
($\delta>0$) under $P^0$, by Yokoyama (1980), we obtain
\begin{equation}\label{eq:uniformbound}
E^0\{ |y_n-E^0(y_n)|^{8+\delta} \} \le K <\infty, ~~~~~ \delta>0
\end{equation}
uniformly in $n$, provided that $\{ \tau_k \}$ is exponential
$\alpha$-mixing, which is proved in Lemma \ref{lem:mixing}.

Therefore, we can apply Theorem \ref{theo:general} to the ACD(1,1)
model and the result follows. $\Box$

\textbf{Proof of Lemma \ref{lem:partialsumLMSD}:}

We present the proof for the case $0 < d < \frac{1}{2}$. The proof
for the case $d=0$ follows along similar lines. Also, we assume here
that $p$ is a positive even integer. The result for all positive odd
integers follows by H\"{o}lder's inequality.

Let $\tilde y_n = y_n - E^0(y_n)$. Since $p \ge 2$ is even and
$E^0{(\tilde y_n)^p}$ can be expressed as a linear combination of
the products of the joint cumulants of $\tilde y_n$ of order $2,
\dots,p$, we have
\begin{eqnarray*}
0 \le E^0|\tilde y_n|^p = E^0(\tilde y_n^p) &=& \sum_\pi \Big [
c_\pi \prod_{j\in\pi} \mbox{cum}( \underbrace{\tilde
y_n,\ldots,\tilde y_n }
_{j ~ \mbox{terms}} ) \Big ] \\
&\le& \sum_\pi \Big [ |c_\pi| \prod_{j\in\pi} |\mbox{cum} (
\underbrace{ \tilde y_n,\ldots,\tilde y_n}_{j ~ \mbox{terms}} ) |
\Big ]
\end{eqnarray*}
where $\pi$ ranges over the additive partitions of $n$ and $c_\pi$
is a finite constant depending on $\pi$.

Since the first order cumulant of $\tilde y_n$ is zero and for all
integers $m \ge 2$, the $m$-th order cumulant of $\tilde y_n$ is
equal to that of $y_n$, it suffices to show that the absolute value
of the $m$-th order cumulant of $y_n$ is bounded uniformly in $n$
under $P^0$, for all $m \in \{ 2,\ldots,p \}$.

We first consider the second and the third order cumulants.

For the second order cumulant $(m=2)$,
\[
|\mbox{cum}(y_n, y_n)| = |\mbox{cum}(\frac{\sum_{k=1}^{n}
\tau_k}{n^{d+\frac{1}{2}}}, \frac{\sum_{s=1}^{n}
\tau_s}{n^{d+\frac{1}{2}}})| \le \frac{1}{n^{2d+1}} \sum_{k=1}^{n}
\sum_{s=1}^{n} |\mbox{cum}(\tau_k, \tau_s)| ~~~~~ .
\]

To calculate the joint cumulant $\mbox{cum}(\tau_k, \tau_s)$, we
briefly introduce some terminology, mainly cited from Brillinger
(1981): consider a (not necessary rectangular) two-way table of
indices,
\begin{eqnarray*}
&(1,1)& ~~~~~ \ldots ~~~~~~~ (1,J_1) \\
&\vdots& ~~~~~ \ldots ~~~~~~~ \vdots \\
&(I,1)& ~~~~~ \ldots ~~~~~~~ (I,J_I)
\end{eqnarray*}
and a partition $P_1 \cup P_2 \cup \ldots \cup P_M$ of its entries.
We say sets $P_{m^\prime}$, $P_{m^{\prime\prime}}$ of the partition
\textbf{hook} if there exist $(i_1,j_1) \in P_{m^\prime}$ and
$(i_2,j_2) \in P_{m^{\prime\prime}}$ such that $i_1=i_2$, i.e. at
least one entry of $P_{m^\prime}$ and one entry of
$P_{m^{\prime\prime}}$ come from the same row in the two-way table.
We say that sets $P_{m^{\prime}}$ and $P_{m^{\prime\prime}}$
\textbf{communicate} if there exists a sequence of sets
$P_{m_1}=P_{m^{\prime}},P_{m_2},\ldots,P_{m_N}=P_{m^{\prime\prime}}$
such that $P_{m_n}$ and $P_{m_{n+1}}$ hook for $n=1,\ldots,N-1$. So
$P_{m^{\prime}}$ and $P_{m^{\prime\prime}}$ communicate as long as
one can find an ordered sequence of sets such that all neighboring
pairs hook, and this sequence links $P_{m^{\prime}}$ and
$P_{m^{\prime\prime}}$ together. Finally a partition is said to be
\textbf{indecomposable} if all sets in the partition communicate.

By Brillinger (1981), Theorem 2.3.2, for a two-way array of random
variables $X_{ij}$, $j=1,\ldots,J_i$, $i=1,\ldots, I$ (see the
corresponding two-way table above), the joint cumulant of the $I$
row products
\[
Y_i = \prod_{j=1}^{J_i} X_{ij}, ~~~~~ i=1,\ldots,I
\]
is given by,
\[
\mbox{cum}(Y_1,\ldots,Y_I) = \sum_\nu \mbox{cum}(X_{ij};ij \in
\nu_1) \ldots \mbox{cum}(X_{ij};ij \in \nu_w)
\]
where the summation is over all indecomposable partition $\nu=\nu_1
\cup \ldots \cup \nu_w$ of the two-way table of indices.

It is more convenient to write the partitions in terms of symbols
representing the random variables, instead of the indices
themselves. We will always use distinct symbols, so that there is a
one-to-one correspondence between the indices and the symbols.
Nevertheless, the random variables represented by distinct symbols
need not be distinct. For example, $e^{h_k}$ and $e^{h_s}$ are
distinct symbols, but if $k=s$, they are not different random
variables. Ultimately, the cumulants are computed from the random
variables.

To compute $\mbox{cum}(\tau_k, \tau_s)$, we use the two-way table of
indices (left) and the corresponding table of symbols (right),
\begin{eqnarray*}
&(1,1)& ~~~~~ (1,2) ~~~~~~~~~~~~~~~~~~~~~~~~~~~~~~~ e^{h_k} ~~~~~ \epsilon_k \\
&(2,1)& ~~~~~ (2,2) ~~~~~~~~~~~~~~~~,~~~~~~~~~~~~~~ e^{h_s} ~~~~~
\epsilon_s
\end{eqnarray*}
with $I=2, J_1=2$ and $J_2=2$.

From Brillinger (1981), Theorem 2.3.1, all joint cumulants
corresponding to partitions with at least one of the symbols
representing $\{ e^{h_k} \}$ and at least one of the symbols
representing $\{ \epsilon_k \}$ in the same set, are zero because
the corresponding random variable sequences are mutually
independent. So for $m=2$, excluding those with at least one of
$e^{h_k},e^{h_s}$ and at least one of $\epsilon_k,\epsilon_s$ in the
same set, the only possible indecomposable partitions (here, the
partition is given in terms of the symbols) are:
\begin{eqnarray*}
&\{e^{h_k}, e^{h_s}\}, \{\epsilon_k, \epsilon_s\}&\\
&\{e^{h_k}, e^{h_s}\}, \{\epsilon_k\}, \{\epsilon_s\}&\\
&\{e^{h_k}\}, \{e^{h_s}\}, \{\epsilon_k, \epsilon_s\}& ~~~~~ .
\end{eqnarray*}

Thus, $|\mbox{cum}(y_n, y_n)| \le A + B + C$, where,
\begin{eqnarray*}
A&=&\frac{1}{n^{2d+1}} \sum_{k=1}^{n} \sum_{s=1}^{n}
|\mbox{cum}(e^{h_k},e^{h_s})| |\mbox{cum}(\epsilon_k,\epsilon_s)| \\
B&=&\frac{1}{n^{2d+1}} \sum_{k=1}^{n} \sum_{s=1}^{n}
|\mbox{cum}(e^{h_k}) \mbox{cum}(e^{h_s})| |\mbox{cum}(\epsilon_k,\epsilon_s)| \\
C&=&\frac{1}{n^{2d+1}} \sum_{k=1}^{n} \sum_{s=1}^{n}
|\mbox{cum}(e^{h_k},e^{h_s})| |\mbox{cum}(\epsilon_k)|
|\mbox{cum}(\epsilon_s)|
\end{eqnarray*}

Both $A$ and $B$ reduce to a single summation because of the
serial independence of the $\{ \epsilon_k \}$, so $A=O(n^{-2d})$
and $B=O(n^{-2d})$. For $C$, by Surgailis and Viano (2002),
Corollary 5.3,
\[
|\mbox{cum}(e^{h_k},e^{h_s})|=e^{\sigma_h^2}|e^{r_{|k-s|}}-1|
\]
where $r_{|k-s|} = \textmd{cov} (h_k,h_s)$ and $\sigma_h^2=
\textmd{var} (h_k)$.

By the assumption on $\{ b_j \}$ in the Theorem \ref{theo:LMSD}, it
follows that $r_s \sim K s^{2d-1}$, as $s \rightarrow \infty$, where
$K >0$, so that
\begin{eqnarray*}
\sum_{k=1}^{n} \sum_{s=1}^{n} |e^{r_{|k-s|}}-1| &\le& 2
\sum_{k=1}^{n} \sum_{s > k}^{n} |e^{r_{|k-s|}}-1| + n |e^{r_{0}}-1|
\\ &\le& K n \sum_{j=1}^{n} j^{2d-1} + n |e^{r_{0}}-1| = O(n^{2d+1}).
\end{eqnarray*}
Thus term $C$ is $O(1)$. Hence, $|\mbox{cum}(y_n, y_n)|$ is $O(1)$.

Next, for the third order cumulant $(m=3)$, we have
\[
|\mbox{cum}(y_n, y_n, y_n)| = \frac{1}{n^{3d+\frac{3}{2}}}
|\sum_{k=1}^{n} \sum_{s=1}^{n} \sum_{u=1}^{n} \mbox{cum}(\tau_k,
\tau_s, \tau_u)| \le \frac{1}{n^{3d+\frac{3}{2}}} \sum_{k=1}^{n}
\sum_{s=1}^{n} \sum_{u=1}^{n} |\mbox{cum}(e^{h_k}\epsilon_k,
e^{h_s}\epsilon_s, e^{h_u}\epsilon_u)| ~~~~~ .
\]

We will use the following two-way table:
\begin{eqnarray*}
e^{h_k} ~~~~~ \epsilon_k \\
e^{h_s} ~~~~~ \epsilon_s \\
e^{h_u} ~~~~~ \epsilon_u
\end{eqnarray*}

For convenience, we group the indecomposable partitions according to
how many sets $(L=1,2,3)$ the symbols $e^{h_k},e^{h_s},e^{h_u}$ are
partitioned into.

We have three groups of indecomposable partitions, excluding those
with at least one of $e^{h_k},e^{h_s},e^{h_u}$ and at least one of
$\epsilon_k,\epsilon_s,\epsilon_u$ in the same set: \\
i) Group 1
\[\{e^{h_k},e^{h_s},e^{h_u}\},\{\epsilon_k,\epsilon_s,\epsilon_u\}\]
\[\{e^{h_k},e^{h_s},e^{h_u}\},\{\epsilon_k,\epsilon_s\},\{\epsilon_u\}\]
\[\{e^{h_k},e^{h_s},e^{h_u}\},\{\epsilon_k,\epsilon_u\},\{\epsilon_s\}\]
\[\{e^{h_k},e^{h_s},e^{h_u}\},\{\epsilon_k\},\{\epsilon_s,\epsilon_u\}\]
\[\{e^{h_k},e^{h_s},e^{h_u}\},\{\epsilon_k\},\{\epsilon_s\},\{\epsilon_u\}\]
ii) Group 2
\[\{e^{h_k},e^{h_s}\},\{e^{h_u}\},\{\epsilon_k,\epsilon_s,\epsilon_u\}\]
\[\{e^{h_k},e^{h_s}\},\{e^{h_u}\},\{\epsilon_k\},\{\epsilon_s,\epsilon_u\}\]
\[\{e^{h_k},e^{h_s}\},\{e^{h_u}\},\{\epsilon_s\},\{\epsilon_k,\epsilon_u\}\]
\[\{e^{h_k},e^{h_u}\},\{e^{h_s}\},\{\epsilon_k,\epsilon_s,\epsilon_u\}\]
\[\{e^{h_k},e^{h_u}\},\{e^{h_s}\},\{\epsilon_k\},\{\epsilon_s,\epsilon_u\}\]
\[\{e^{h_k},e^{h_u}\},\{e^{h_s}\},\{\epsilon_s\},\{\epsilon_k,\epsilon_u\}\]
\[\{e^{h_u},e^{h_s}\},\{e^{h_k}\},\{\epsilon_k,\epsilon_s,\epsilon_u\}\]
\[\{e^{h_u},e^{h_s}\},\{e^{h_k}\},\{\epsilon_k\},\{\epsilon_s,\epsilon_u\}\]
\[\{e^{h_u},e^{h_s}\},\{e^{h_k}\},\{\epsilon_s\},\{\epsilon_k,\epsilon_u\}\]
iii) Group 3
\[\{e^{h_k}\},\{e^{h_s}\},\{e^{h_u}\},\{\epsilon_k,\epsilon_s,\epsilon_u\}.\]

We next study the order of the dominant contribution to
$|\mbox{cum}(y_n,y_n,y_n)|$ corresponding to each group.

In Group 1, the dominant term arises from the last partition since
it yields a triple summation,
\[
\frac{1}{n^{3d+\frac{3}{2}}} \sum_{k=1}^{n} \sum_{s=1}^{n}
\sum_{u=1}^{n}
|\mbox{cum}(e^{h_k},e^{h_s},e^{h_u})||\mbox{cum}(\epsilon_k)||\mbox{cum}(\epsilon_s)||\mbox{cum}(\epsilon_u)|
= \frac{\mu_{\epsilon}^3}{n^{3d+\frac{3}{2}}} \sum_{k=1}^{n}
\sum_{s=1}^{n} \sum_{u=1}^{n} |\mbox{cum}(e^{h_k},e^{h_s},e^{h_u})|
\]
where $\mu_\epsilon=E^0(\epsilon_1)$.

By Surgailis and Viano (2002), Corollary 5.3,
\begin{eqnarray*}
&& \sum_{k=1}^{n} \sum_{s=1}^{n} \sum_{u=1}^{n}
|\mbox{cum}(e^{h_k},e^{h_s},e^{h_u})| \\ &\le& \sum_{k=1}^{n}
\sum_{s=1}^{n} \sum_{u=1}^{n} e^{\frac{3}{2}\sigma_h^2}
|e^{r_{|k-s|}}-1||e^{r_{|k-u|}}-1||e^{r_{|s-u|}}-1|
\\ &+& \sum_{k=1}^{n}
\sum_{s=1}^{n} \sum_{u=1}^{n} e^{\frac{3}{2}\sigma_h^2}
|e^{r_{|k-s|}}-1||e^{r_{|k-u|}}-1| + \sum_{k=1}^{n} \sum_{s=1}^{n}
\sum_{u=1}^{n} e^{\frac{3}{2}\sigma_h^2}
|e^{r_{|k-s|}}-1||e^{r_{|s-u|}}-1|
\\ &+& \sum_{k=1}^{n} \sum_{s=1}^{n} \sum_{u=1}^{n} e^{\frac{3}{2}\sigma_h^2}
|e^{r_{|k-u|}}-1||e^{r_{|s-u|}}-1|
\end{eqnarray*}

The last three summations are actually the same due to symmetry: we
can simply relabel the indices in the last summation by $s
\leftrightarrow u$. As for the first summation, since $|r_{|k-u|}|=|
\textmd{cov} (h_k,h_u)| \le \sigma_h^2 = \textmd{var} (h_k)$, we
have $|e^{r_{|k-u|}}-1| \le (e^{\sigma_h^2}+1) < \infty$. So
\begin{eqnarray*}
\sum_{k=1}^{n} \sum_{s=1}^{n} \sum_{u=1}^{n}
|\mbox{cum}(e^{h_k},e^{h_s},e^{h_u})| &\le& K \sum_{k=1}^{n}
\sum_{s=1}^{n} \sum_{u=1}^{n} e^{\frac{3}{2}\sigma_h^2}
|e^{r_{|k-s|}}-1||e^{r_{|s-u|}}-1|
\\ &+& 3 \sum_{k=1}^{n} \sum_{s=1}^{n}
\sum_{u=1}^{n} e^{\frac{3}{2}\sigma_h^2}
|e^{r_{|k-s|}}-1||e^{r_{|s-u|}}-1| \\ &\le& K \sum_{k=1}^{n}
\sum_{s=1}^{n} \sum_{u=1}^{n} |e^{r_{|k-s|}}-1||e^{r_{|s-u|}}-1|
~~~~~ (\mbox{for~some}~K > 0)
\\ &=& O(n^{4d+1})
\end{eqnarray*}

The last step follows from Lemma \ref{lem:tree}. So
$\frac{\mu_{\epsilon}^3}{n^{3d+\frac{3}{2}}} \sum_{k=1}^{n}
\sum_{s=1}^{n} \sum_{u=1}^{n} |\mbox{cum}(e^{h_k},e^{h_s},e^{h_u})|$
converges to zero because $(4d+1)<(3d+\frac{3}{2})$.

Similarly, the dominant contribution from Group 2 is of order
\[
\frac{1}{n^{3d+\frac{3}{2}}} \sum_i \sum_j
|\mbox{cum}(e^{h_i},e^{h_j})| |\mbox{cum}(e^{h_{j}})|
\]
Note that in Group 2, all three of $e^{h_k},e^{h_s},e^{h_u}$ are
partitioned into two sets. Therefore, partitions with all three of
$\epsilon_k,\epsilon_s,\epsilon_u$ in different sets are not
indecomposable, so the dominant contribution is a double sum,
\[
\frac{1}{n^{3d+\frac{3}{2}}} \sum_i \sum_j
|\mbox{cum}(e^{h_i},e^{h_j})||\mbox{cum}(e^{h_{j}})| =
\frac{\mu_{e^h}}{n^{3d+\frac{3}{2}}} \sum_i \sum_j
|\mbox{cum}(e^{h_i},e^{h_j})| \le K n^{(2d+1)-(3d+\frac{3}{2})} =
O(n^{-d-\frac{1}{2}})
\]
where $\mu_{e^h} = E^0(e^{h_1})$.

So the dominant term in Group 2 also converges to zero.

For Group 3, all three of $e^{h_k},e^{h_s},e^{h_u}$ are partitioned
into three different sets, so that the part of the partition
involving $\epsilon_k,\epsilon_s,\epsilon_u$ must be
$\{\epsilon_k,\epsilon_s,\epsilon_u\}$ in order to be
indecomposable. The resulting summation now is only a single one of
order $O(n^1)$. The dominant contribution again converges to zero.

Notice that the order of the dominant contribution from group 3
($O(n^{-3d-\frac{1}{2}})$) is of smaller order than that from group
2 ($O(n^{-d-\frac{1}{2}})$), which is of smaller order of that from
group 1 ($O(n^{d-\frac{1}{2}})$). This will be shown to hold in
general for any $m$-th order joint cumulant.

Next, we prove that the $m$-th order joint cumulant, which satisfies
\begin{equation}\label{eq:cum}
| \mbox{cum}(\underbrace{y_n,\ldots,y_n}_{m~\mbox{terms}}) | \le
\frac{1}{n^{m(d+\frac{1}{2})}} \sum_{k_1=1}^n \ldots
\sum_{k_m=1}^n|\mbox{cum}(e^{h_{k_1}} \epsilon_{k_1}, \ldots,
e^{h_{k_m}} \epsilon_{k_m})|
\end{equation}
converges to zero for all $m>2$.

The indecomposable partitions of $(e^{h_{k_1}}\epsilon_{k_1},
\ldots, e^{h_{k_m}}\epsilon_{k_m})$ are organized in a similar
manner as before into $m$ groups, where in Group $L$ the symbols
$e^{h_{k_1}}, \ldots, e^{h_{k_m}}$ are divided into $L$ sets
($L=1,\dots,m$).

a) First, consider Group 1. The dominant contribution to the
righthand side of (\ref{eq:cum}) corresponding to Group 1 must be
the one from the partition in which all of the symbols
$e^{h_{k_1}},\ldots,e^{h_{k_m}}$ are in one set and each of the
symbols $\epsilon_{k_1}, \ldots, \epsilon_{k_m}$ is in a set by
itself. The resulting summation is an $m$-fold summation. By
Corollary 5.3 of Surgailis and Viano (2002), the absolute value of
the $m$-th joint cumulant,
$|\mbox{cum}(e^{h_{k_1}},\ldots,e^{h_{k_m}})|$, is bounded by a
summation taken over all connected graphs with $m$ vertices. Each
entry of the summation is a product of terms of the form
$|e^{r_{|k_i-k_j|}}-1|$ along the edges that connect vertices
$k_i$ and $k_j$ of a connected $m$-vertex graph.

For a graph with $m$ vertices, we need at least $(m-1)$ edges to
connect them. It is known (see Andrasfai, 1977, Chapter 2) that any
connected $m$-vertex graph with $(m-1)$ edges may be represented as
a \textit{tree}. Let $W_{\{k_i,\ldots,k_j\}} < \infty$ be the total
number of trees with vertices labeled by $k_i, k_{i+1}, \ldots,
k_j$.

If a connected $m$-vertex graph used in applying Corollary 5.3 of
Surgailis and Viano (2002) has more than $(m-1)$ edges, it is not
a tree, and there will be more than $(m-1)$ terms of the form
$|e^{r_{|k_i-k_j|}}-1|$ being multiplied together in the $m$-fold
summation in (\ref{eq:cum}). But, for all $k_i,k_j$,
$|r_{|k_i-k_j|}|=| \textmd{cov} (h_{k_i},h_{k_j})| \le \sigma_h^2
= \textmd{var} (h_{k_i})$, so $|e^{r_{|k_i-k_j|}}-1| \le
(e^{\sigma_h^2}+1) < \infty$, and for any connected $m$-vertex
graph with more than $(m-1)$ edges, there exists an $m$-vertex
subgraph that has a tree representation. So we can retain a
product of $(m-1)$ terms of the form $|e^{r_{|k_i-k_j|}}-1|$ in
the $m$-fold summation in (\ref{eq:cum}) and move remaining terms
out of the summation, bounding each by $(e^{\sigma_h^2}+1)$. The
resulting product of $(m-1)$ terms of the form
$|e^{r_{|k_i-k_j|}}-1|$ is itself a product over the edges of an
$m$-vertex tree.

In all, $|\mbox{cum}(e^{h_{k_1}},\ldots,e^{h_{k_m}})|$ is bounded by
a constant times a summation over the set $G_{\{k_1,\ldots,k_m\}}$
of all $W_{\{k_1,\ldots,k_m\}}$ trees. Each entry of the summation
is a product of terms of the form $|e^{r_{|k_i-k_j|}}-1|$ being
multiplied over the $(m-1)$ edges of the tree. Thus, we have
\begin{eqnarray*}
\sum_{k_1=1}^n \ldots \sum_{k_m=1}^n | \mbox{cum}(e^{h_{k_1}},
\ldots, e^{h_{k_m}})| &\le& K \sum_{k_1=1}^n \ldots \sum_{k_m=1}^n
\Big \{ \sum_{G_{\{k_1,\ldots,k_m\}}} ~ \prod_{(k_i,k_j) \in
\Omega(G_{\{k_1,\ldots,k_m\}})} |e^{r_{|k_i-k_j|}}-1| \Big \} , ~~~
(K
>0) \\ &=& K \sum_{G_{\{k_1,\ldots,k_m\}}}
\sum_{k_1=1}^n \ldots \sum_{k_m=1}^n \Big \{
\underbrace{\prod_{(k_i,k_j) \in \Omega(G_{\{k_1,\ldots,k_m\}})}
|e^{r_{|k_i-k_j|}}-1|}_{(m-1)~\mbox{terms}} \Big \}
\end{eqnarray*}
where $\Omega(G_{\{k_1,\ldots,k_m\}})$ is the set of edges of the
graph indexed by $G_{\{k_1,\ldots,k_m\}}$.

By Lemma \ref{lem:tree}, each entry of the summation over
$G_{\{k_1,\ldots,k_m\}}$ is of order $O(n^{2dm-2d+1})$. Also this
summation is taken over a finite number of graphs
($W_{\{k_1,\ldots,k_m\}} < \infty$), therefore
\[
\sum_{k_1=1}^n \ldots \sum_{k_m=1}^n
|\mbox{cum}(e^{h_{k_1}},\ldots,e^{h_{k_m}})| = O(n^{2dm-2d+1})_\cdot
\]

Because the normalization term in (\ref{eq:cum}) is of order
$O(n^{m(d+\frac{1}{2})})$, the dominant contribution to
$\mbox{cum}(\underbrace{y_n, \ldots, y_n}_{m ~ \mbox{terms}})$ from
Group 1 converges to zero, for any $m>2$.

b) For Group 2, the symbols $e^{h_{k_1}},\ldots,e^{h_{k_m}}$ are
partitioned into two sets. Thus, the partitions with each of the $m$
symbols $\epsilon_{k_1}, \ldots, \epsilon_{k_m}$ in a set by itself
are not indecomposable. Relabel the two sets as $\{ e^{h_{g_1}},
\ldots, e^{h_{g_q}} \}$, $\{ e^{h_{g_{q+1}}}, \ldots, e^{h_{g_m}}
\}$. Since the partition must be indecomposable, there must be one
$I \in (1,\ldots,q)$ and one $J \in (q+1,\ldots,m)$, such that $g_I
=g_J$. The dominant contribution to (\ref{eq:cum}) from Group 2 is
therefore
\begin{equation}\label{eq:cumg2}
\frac{1}{n^{m(d+\frac{1}{2})}} \sum_{g_1=1}^n \ldots \sum_{g_m=1}^n
|\mbox{cum}(e^{h_{g_1}},\ldots,e^{h_{g_q}})|
|\mbox{cum}(e^{h_{g_{q+1}}},\ldots,e^{h_{g_m}})|
|\mbox{cum}(\epsilon_{g_I},\epsilon_{g_J})|
\end{equation}

Similarly as above, after applying Corollary 5.3 of Surgailis and
Viano (2002) and after bounding certain terms, we obtain

\begin{eqnarray*}
&&\sum_{g_1=1}^n \ldots \sum_{g_m=1}^n
|\mbox{cum}(e^{h_{g_1}},\ldots,e^{h_{g_q}})|
|\mbox{cum}(e^{h_{g_{q+1}}},\ldots,e^{h_{g_m}})|
|\mbox{cum}(\epsilon_{g_I},\epsilon_{g_J})|
\\
&\le& K \sum_{g_1=1}^n \ldots \sum_{g_m=1}^n \Big \{
\sum_{G_{\{g_1,\ldots,g_q\}}} ~ \underbrace{\prod_{(g_i,g_j) \in
\Omega(G_{\{g_1,\ldots,g_q\}})}
|e^{r_{|g_i-g_j|}}-1|}_{(q-1)~\mbox{terms}} \Big\} \\ &&
~~~~~~~~~~~~~~~~~~ \cdot \Big\{ \sum_{G_{\{g_{q+1},\ldots,g_m\}}} ~
\underbrace{ \prod_{(g_i,g_j) \in
\Omega(G_{\{g_{q+1},\ldots,g_m\}})}
|e^{r_{|g_i-g_j|}}-1|}_{(m-q-1)~\mbox{terms}} \Big\} \Big\{
|\mbox{cum}(\epsilon_{g_I},\epsilon_{g_J})| \Big \} \\
&=& K \sum_{G_{\{g_1,\ldots,g_q\}}}
\sum_{G_{\{g_{q+1},\ldots,g_m\}}} \sum_{g_1=1}^n \ldots
\sum_{g_m=1}^n \textbf{1}_{\{g_I=g_J\}} \\ &&
~~~~~~~~~~~~~~~~~~~~~~~~~~~~~~~~ \cdot \Big \{ \underbrace{
\prod_{(g_i,g_j) \in \Omega(G_{\{g_1,\ldots,g_q\}})}
|e^{r_{|g_i-g_j|}}-1|  \prod_{(g_i,g_j) \in
\Omega(G_{\{g_{q+1},\ldots,g_m\}})}
|e^{r_{|g_i-g_j|}}-1|}_{(m-2)~\mbox{terms, denote
as}~\Gamma(g_1,\ldots,g_m:
G_{\{g_1,\ldots,g_q\}},G_{\{g_{q+1},\ldots,g_m\}})} \Big \}_.
\end{eqnarray*}

As mentioned before, any graph $G_a$ in $G_{\{g_1,\ldots,g_q\}}$ and
any graph $G_b$ in $G_{\{g_{q+1},\ldots,g_m\}}$, can be represented
by trees with $q$ and $(m-q)$ vertices, respectively. Since for any
two trees, the resulting structure obtained by merging one vertex
from each tree is again a tree, under the constraint $g_I=g_J$,
there exists a graph $G_c$ in
$G_{\{g_1,\ldots,g_{I-1},g_{I+1},\ldots,g_m\}}$, such that $G_c$ is
obtained by merging $G_a$ and $G_b$ together at the vertex
$g_I=g_J$.

Therefore, the numerical value of the term $\Gamma$ evaluated for
graphs $G_a$ and $G_b$ and indices $\{ g_1,\ldots,g_m\}$ with the
constraint $g_I=g_J$ (which follows from the independence of the $\{
\epsilon_{g_i} \}$) is equal to the value of the term $\Phi$
(defined below) evaluated using the graph $G_c$ in
$G_{\{g_1,\ldots,g_{I-1},g_{I+1},\ldots,g_m\}}$ and indices $\{
g_1,\ldots,g_{I-1},g_{I+1},\ldots,g_m \}$ without any constraint on
the values of these indices. After re-parameterizing $\{
g_1,\ldots,g_{I-1},g_{I+1},\ldots,g_m \}$ by
$\{l_1,\ldots,l_{m-1}\}$, we obtain
\begin{eqnarray*}
&&\sum_{g_1=1}^n \ldots \sum_{g_m=1}^n
|\mbox{cum}(e^{h_{g_1}},\ldots,e^{h_{g_q}})|
|\mbox{cum}(e^{h_{g_{q+1}}},\ldots,e^{h_{g_m}})|
|\mbox{cum}(\epsilon_{g_I},\epsilon_{g_J})|
\\ &\le& K \sum_{G_{\{l_1,\ldots,l_{m-1}\}}} \sum_{l_1=1}^n \ldots
\sum_{l_{m-1}=1}^n \underbrace{ \prod_{(l_i,l_j) \in
\Omega(G_{\{l_1,\ldots,l_{m-1}\}})}
|e^{r_{|l_i-l_j|}}-1|}_{(m-2)~\mbox{terms, denote
as}~\Phi(l_1,\ldots,l_{m-1}: G_{\{l_1,\ldots,l_{m-1}\}})} \\ &=&
O(n^{2d(m-2)+1})
\end{eqnarray*}
where the final equality follows from Lemma \ref{lem:tree}.

The above $(m-1)$-fold summation for Group 2 is of smaller order
than the $m$-fold summation from Group 1, which was
$O(n^{2d(m-1)+1})$. Hence, the dominant contribution from Group 2
also converges to zero.

c) In general, for Group $L \in \{ 1,\ldots,m \}$, the symbols
$e^{h_{k_1}},\ldots,e^{h_{k_m}}$ are partitioned into $L$ sets.
Relabel the $L$ sets as $\{ e^{h_{g_1}}, \ldots, e^{h_{g_{q_1}}}
\}$, $\{ e^{h_{g_{q_1+1}}}, \ldots, e^{h_{g_{q_2}}} \}$, \ldots, $\{
e^{h_{g_{q_{L-1}+1}}}, \ldots, e^{h_{g_m}} \}$. Since the partition
must be indecomposable, there must be $L$ indices $\{ I,J,\ldots,Z
\}$, where $I \in (1,\ldots,q_1), J \in (q_1+1,\ldots, q_2), \ldots,
Z \in (q_{L-1}+1,\ldots,m)$, such that $\underbrace{g_I=g_J=\ldots =
g_Z}_{L~\mbox{terms}}$. The dominant contribution to (\ref{eq:cum})
from Group $L$ is then,
\begin{equation}\label{eq:cumL}
\frac{1}{n^{m(d+\frac{1}{2})}} \sum_{g_1=1}^n \ldots \sum_{g_m=1}^n
\underbrace{|\mbox{cum}(e^{h_{g_1}},\ldots,e^{h_{g_{q_1}}})| \ldots
|\mbox{cum}(e^{h_{g_{q_{L-1}+1}}},\ldots,e^{h_{g_m}})|}_{L-\mbox{terms}}
|\mbox{cum}(\underbrace{\epsilon_{g_I},\epsilon_{g_J},\ldots,\epsilon_{g_Z}}_{L~\mbox{terms}})|_.
\end{equation}

Similarly as before, we obtain
\begin{eqnarray*}
&& \sum_{g_1=1}^n \ldots \sum_{g_m=1}^n
\underbrace{|\mbox{cum}(e^{h_{g_1}},\ldots,e^{h_{g_{q_1}}})| \ldots
|\mbox{cum}(e^{h_{g_{q_{L-1}+1}}},\ldots,e^{h_{g_m}})|}_{L-\mbox{terms}}
|\mbox{cum}(\underbrace{\epsilon_{g_I},\epsilon_{g_J},\ldots,\epsilon_{g_Z}}_{L~\mbox{terms}})|
\\
&\le& K \underbrace{ \sum_{G_{\{g_1,\ldots,g_{q_1}\}}} \ldots
\sum_{G_{\{g_{q_{L-1}+1},\ldots,g_m\}}} }_{L-\mbox{fold}}
\sum_{g_1=1}^n \ldots \sum_{g_m=1}^n
\textbf{1}_{\{\underbrace{g_I=g_J=\ldots=g_Z}_{L~\mbox{terms}}\}} \\
&& \cdot \underbrace{ \Big\{ \prod_{(g_i,g_j) \in
\Omega(G_{\{g_1,\ldots,g_{q_1}\}})} |e^{r_{|g_i-g_j|}}-1| ~~~ \ldots
~~~ \prod_{(g_i,g_j) \in \Omega(G_{\{g_{q_{L-1}+1},\ldots,g_m\}})}
|e^{r_{|g_i-g_j|}}-1|\Big\} }_{(m-L) ~ \mbox{terms}}  \\ &\le& K
\sum_{G_{\{l_1,\ldots,l_{m-L+1}\}}} \sum_{l_1=1}^n \ldots
\sum_{l_{m-L+1}=1}^n \underbrace{ \prod_{(l_i,l_j) \in
\Omega(G_{\{l_1,\ldots,l_{m-L+1}\}})}
|e^{r_{|l_i-l_j|}}-1|}_{(m-L)~\mbox{terms}} \\ &=& O(n^{2d(m-L)+1}),
\end{eqnarray*}
by Lemma \ref{lem:tree}.

The constraint $\underbrace{g_I=g_J = \ldots=g_Z}_{L~\mbox{terms}}$
allows the re-parameterization from $\{ g_1,\ldots,g_m \}$ to $\{
l_1,\ldots,l_{m-L+1} \}$ and reduces the $m$-fold summation in
(\ref{eq:cumL}) to an $(m-L+1)$-fold summation in the last
inequality. It was shown for Group 2 that the graph obtained by
merging one vertex from each of any pair of trees is again a tree.
By induction, we obtain a tree by merging one vertex from each of
$L>2$ trees, which allows us to apply Lemma \ref{lem:tree} with
$M=m-L+1$ in the last step.

So, the dominant contribution from Group $L$ is
$O(n^{2d(m-L)+1-m(d+\frac{1}{2})})$, ($L=1,\ldots,m$). Since $d>0$,
the dominant contribution from \textit{all} groups occurs for $L=1$.
Finally, the dominant contribution from Group 1 is
$O(n^{2d(m-1)+1-m(d+\frac{1}{2})})$, which tends to zero for $m>2$
since $d<\frac{1}{2}$. $\Box$

\begin{lem} \label{lem:supz}
For durations $\{ \tau_k \}$ satisfying the assumptions of Theorem
\ref{theo:general},
\[
\limsup_{t} E[Z^4(t)] < \infty
\]
where $Z(t)$ is defined by Equation (\ref{eq:defZt}).
\end{lem}

\textbf{Proof:} By Chung (1974, Theorem 3.2.1, page 42), $E[Z(t)^4]
\le 1 + \sum_{s=1}^\infty P[Z^4(t) \ge s]$. Thus, it suffices to
show that
\begin{equation}\label{eq:sufficient}
\limsup_t \sum_{s=1}^\infty P[Z^4(t) \ge s] < \infty.
\end{equation}

Note that for any real $k$,
\begin{equation} \label{eq:ConditionA}
N(t) \ge k \Longleftrightarrow \sum_{i=1}^{ \lfloor k \rfloor} u_i
\le t.
\end{equation}

We have
\begin{equation}
P[Z^4(t) \ge s] = P[Z(t) \le -s^{1/4}] + P[Z(t) \ge s^{1/4}].
\end{equation}
Consider the second term $P[Z(t) \ge s^{1/4}]$. Using
(\ref{eq:ConditionA}), we obtain
\begin{eqnarray*}
P[Z(t) \ge s^{1/4}] &=& P[N(t) \ge \frac{t}{\mu} + s^{1/4}t^{1/2+d}]
\\ &=& P(\sum_{i=1}^{\lfloor g(t,s) \rfloor} u_i \le t)
\end{eqnarray*}
where $g(t,s)=\frac{t}{\mu}+s^{1/4}t^{1/2+d}$.

So,
\begin{eqnarray*}
P[Z(t) \ge s^{1/4}] &=& P(\sum_{i=1}^{\lfloor g(t,s) \rfloor} u_i
\le t) \\ &=& P\Big(\frac{\sum_{i=1}^{\lfloor g(t,s) \rfloor}
(u_i-\mu)}{\lfloor g(t,s) \rfloor^{1/2+d}}
\le \frac{t-\lfloor g(t,s) \rfloor\mu}{\lfloor g(t,s) \rfloor^{1/2+d}}\Big). \\
\end{eqnarray*}

Denote
\[
U=\frac{\sum_{i=1}^{\lfloor g(t,s) \rfloor} (u_i-\mu)}{\lfloor
g(t,s) \rfloor^{1/2+d}} \twoem .
\]
Since $-\lfloor x \rfloor <-x+1$ for $x>0$, we obtain for any
positive $p$,
\begin{eqnarray*}
P[Z(t) \ge s^{1/4}] &=& P\Big(U \le \frac{t-\lfloor g(t,s)
\rfloor\mu}{\lfloor g(t,s) \rfloor^{1/2+d}}\Big) \\ &\le& P\Big(U
\le \frac{-\mu s^{1/4}t^{1/2+d} + \mu}{\lfloor g(t,s) \rfloor^{1/2+d}}\Big) \\
&\le& P\Big(|U| \ge \frac{\mu s^{1/4}t^{1/2+d} - \mu}{\lfloor g(t,s)
\rfloor^{1/2+d}}\Big) \\ &\le& K
E(|U|^p)\frac{(\frac{t}{\mu}+s^{1/4}t^{1/2+d})^{(1/2+d)p}}{\mu^p[
s^{1/4}t^{1/2+d}-1]^p}
\end{eqnarray*}

For $t \ge 4$, since $s^{1/4}t^{1/2+d}-1 \ge
\frac{1}{2}s^{1/4}t^{1/2+d}$ and $\frac{1}{2}+d < 1 $, we obtain
\begin{equation}\label{eq:greater}
P[Z(t) \ge s^{1/4}] \le K E(|U|^p)
\frac{1}{s^{\frac{p}{4}(\frac{1}{2}-d)}}
\end{equation}

Now, consider
\[
P[Z(t) \le -s^{1/4}] = P[N(t) \le \frac{t}{\mu} - s^{1/4}t^{1/2+d}]
\twoem.
\]

Let $a(t)=\frac{t^{2-4d}}{\mu^4}$ and
$v(t,s)=\frac{t}{\mu}-s^{1/4}t^{1/2+d}$. Using
(\ref{eq:ConditionA}), we have
\begin{eqnarray*}
P[Z(t) \le -s^{1/4}] &=& P(\sum_{i=1}^{\lfloor v(t,s) \rfloor} u_i
> t),~~~~~s < a(t)
\\ &=& P(u_1 > t), ~~~~~~~~~~~~~s=a(t) \\ &=&
0,~~~~~~~~~~~~~~~~~~~~~~~~s>a(t)
\end{eqnarray*}

For $s<a(t)$, we have $v(t,s) > 0$. Let
\[
W=\frac{\sum_{i=1}^{\lfloor v(t,s) \rfloor} (u_i-\mu) }{\lfloor
v(t,s) \rfloor^{1/2+d}} \twoem.
\]
Then for any positive $p$,
\begin{eqnarray*}
P[Z(t) \le -s^{1/4}] = P[N(t) \le v(t,s)] &=&
P\Big(\frac{\sum_{i=1}^{\lfloor v(t,s) \rfloor} (u_i-\mu) }{\lfloor
v(t,s) \rfloor^{1/2+d}}
> \frac{t-\mu \lfloor v(t,s) \rfloor}{\lfloor v(t,s) \rfloor^{1/2+d}}\Big) \\ &\le& P\Big(W >
\frac{\mu
s^{1/4}t^{1/2+d}}{(\frac{t}{\mu}-s^{1/4}t^{1/2+d})^{1/2+d}} \Big)
\\ &\le& K E(|W|^p) \frac{(\frac{t}{\mu}-s^{1/4}t^{1/2+d})^{(1/2+d)p}}{s^{p/4}t^{(1/2+d)p}}
\\ &\le& K E(|W|^p) \frac{1}{s^{p/4}}
\end{eqnarray*}
i.e.
\begin{equation}\label{eq:less}
P[Z(t) \le -s^{1/4}] \le K E(|W|^p) \frac{1}{s^{p/4}} \twoem.
\end{equation}

For $s=a(t)$, $P[Z(t) \le -s^{1/4}] = P[u_1>t] \le
\frac{E(u_1)}{t}$.

For $s>a(t)$, $P[Z(t) \le -s^{1/4}] = 0$.

Select any positive $p$ such that $\frac{p}{4}(\frac{1}{2}-d)
> 1$ and thus $\frac{p}{4}>1$ since $0 < \frac{1}{2}-d <1$.
If it can be shown that $\sup_{t>1,s>1} E(|U|^p) < \infty$ and
$\sup_{t>1,s>1} E(|W|^p) < \infty$, then by (\ref{eq:greater}) and
(\ref{eq:less}), it follows that $P[Z^4(t) \ge s]$ is summable,
uniformly in $t$. Thus, (\ref{eq:sufficient}) follows and the proof
is complete.

We next show that indeed $\sup_{t>1,s>1} E(|U|^p) < \infty$ and
$\sup_{t>1,s>1} E(|W|^p) < \infty$ for all positive $p$ when $d \in
(0,\frac{1}{2})$ and for $p=8+\delta,\delta>0$ when $d=0$. Define
\[
B_1 = \frac {u_1 - \mu} {\lfloor g(t,s) \rfloor ^{1/2+d}} \twoem ,
\twoem B_2 = \frac {\sum_{i=2}^{\lfloor g(t,s) \rfloor} (\tau_i -
\mu)} {\lfloor g(t,s) \rfloor ^{1/2+d}},
\]
so that $U=B_1+B_2$. By Minkowski's Inequality,
\[
E[|U|^p] \leq \left [  \left ( E|B_1|^p \right )^{1/p} + \left (
E|B_2|^p \right )^{1/p}\right ] ^p \twoem.
\]
Since $u_1 \leq \tau_1$, using $h(x) = (x+\mu)^p$ in
(\ref{eq:ETau1}), and since by assumption $iv$), $\tau_1$ has all
finite moments up to order $p$ under $P^0$, we have
\[
\sup_{t>1,s>1} E|B_1|^p < \infty \twoem.
\]
From Baccelli and Br$\acute{\text{e}}$maud (2003, Equation 1.2.25,
page 20) that for any measurable function $h$,
\[
E[h(\tau_2 , \ldots , \tau_n)] = \lambda E^0[\tau_1 h(\tau_2
,\ldots,\tau_n )] \twoem.
\]
This, together with the Cauchy-Schwarz inequality, yields
\[
E |B_2|^p = \lambda E^0 (\tau_1 |B_2|^p) \leq \lambda
[E^0(\tau_1^2)]^{1/2} \, [ E^0|B_2|^{2p}\,]^{1/2} \twoem,
\]
where $\lambda = 1/E^0(\tau_1)$. By assumption $iv)$,
$\sup_{t>1,s>1} E^0|B_2|^p <\infty$. for all positive $p$ when $d
\in (0,\frac{1}{2})$ and for $p=8+\delta,\delta>0$ when $d=0$. It
follows that $\sup_{t>1,s>1} E [ |U|^p] < \infty$. By a similar
argument, $\sup_{t>1,s>1} E [ |W|^p] < \infty$. $\Box$

\begin{lem}\label{lem:tree}
For any $M>2$ and $0<d<\frac{1}{2}$,
\begin{equation}\label{eq:lemma2sum}
\underbrace{ \sum_{k_1=1}^n \ldots \sum_{k_M=1}^n }_{M-\mbox{fold}}
\Big \{ \underbrace{ \prod_{(k_i,k_j) \in \Omega(G)}
|e^{r_{|k_i-k_j|}}-1| }_{(M-1)~\mbox{terms}} \Big \} =
O(n^{2d(M-1)+1})
\end{equation}
where $\Omega(G)$ is the set of edges of $G$, $G$ is any connected
$M$-vertex graph with vertices $\{ k_1, \ldots, k_M \}$ and $(M-1)$
edges; $r_{|k_i-k_j|}=\textmd{cov}(h_{k_i}, h_{k_j}), 1 \le i \le M,
1 \le j \le M$, $\{ h_{k_i} \}$ is a long memory process with memory
parameter $d$.
\end{lem}

\textbf{Proof:} Since $G$ is a connected graph with $M$ vertices
and $(M-1)$ edges, it can be represented as a tree (see Andrasfai
1977, Chapter 2). The tree representation is not unique. Fix a
particular representation. Then there is one vertex with no
parent, called the root. A vertex with both a parent and a child
is called a node. A vertex with no child is called a leaf.

We proceed iteratively. First, select any leaf vertex.  By
definition of a leaf, the corresponding index only appears once in
the product, so the sum on this index can be evaluated for this term
only, holding the other terms fixed. Since $r_s \sim C s^{2d-1}$ as
$s \rightarrow \infty$, we have for any fixed integer $i$ with $1
\le i \le n$, $\sum_{j=1}^n |e^{r_{|i-j|}}-1| = O(n^{2d})$.

It follows that the sum on the first index is $O(n^{2d})$. Next,
delete the leaf just used from the tree. The resulting graph is
again a tree. Repeat the process of selecting a leaf, performing the
corresponding sum and deleting the leaf until only the root remains.
The $M$-fold sum in (\ref{eq:lemma2sum}) is now bounded by a
constant times the sum of $n$ terms each of which is
$O(n^{2d(M-1)})$. Thus, the sum in (\ref{eq:lemma2sum}) is
$O(n^{2d(M-1)+1})$. $\Box$

\begin{lem}\label{lem:mixing}
Under the LMSD model described in Theorem \ref{theo:LMSD} with
memory parameter $d \in [0,\frac{1}{2})$, $P^0$ is $\{\tau_k\}$-
mixing; The durations $\{ \tau_k \}$ generated by the ACD(1,1)
model described in Theorem \ref{theo:ACD} are exponential
$\alpha$-mixing.
\end{lem}

\textbf{Proof:} Under $P^0$, $\{h_k\}$ is a stationary Gaussian
process with a log spectral density having an integral on
$[-\pi,\pi]$ that is greater than $-\infty$, so that the
innovation variance is positive. Since Gaussian processes are time
reversible, it follows that we can represent $h_k =
\sum_{j=0}^{\infty} a_j w_{k+j}$ where $\sum a_j^2 < \infty$ and
$\{w_k\}$ is an $iid$ Gaussian sequence. Arguing as in the proof
of Theorem 17.3.1 of Ibragimov and Linnik (1971), $pp$. 311--312,
replacing $\{\ldots w_{k-1}, w_k\}$ by $\{w_k, w_{k+1}, \ldots
\}$, it follows that $P^0$ is $\{h_k\}$-mixing. Since the
$\{\epsilon_k\}$ are $iid$ it follows that $P^0$ is also
$\{\epsilon_k\}$-mixing. Since for any process $\{\xi_k\}$, $P^0$
is $\{\xi_k\}$-mixing if and only if the future tail
$\sigma$-field of $\{\xi_k\}$ is trivial (see, e.g., Nieuwenhuis
(1989), Equation (3.3)), it follows from Lemma \ref{lem:Product}
that $P^0$ is $\{\tau_k\}$-mixing, where $\tau_k = e^{h_k}
\epsilon_k$.

For the ACD(1,1) model, by Proposition 17 of Carrasco and Chen
(2002), $\{ \tau_k \}$ is exponential $\beta$-mixing (or also
called absolutely regular) if $\{ \tau_0, \psi_0 \}$ are
initialized from the stationary distribution. Their result still
holds for a doubly infinite sequence $\{ \tau_k \}, k \in
(-\infty,\infty)$. It is well known that $\beta$-mixing implies
$\alpha$-mixing (or strong mixing), (see Bradley (2005), Section
2.1). Therefore, $\{ \tau_k \}$ is also exponential
$\alpha$-mixing, which further implies $\{ \tau_k \}$-mixing of
$P^0$ for the ACD(1,1) model, see Nieuwenhuis (1989), Equation
(3.5). $\Box$

\begin{lem}\label{lem:Product}
  Let $\{\xi_s\}$ and $\{\zeta_s\}$ be two \emph{independent} processes
  whose future tail $\sigma$-fields are trivial. Then the future
  tail $\sigma$-field of the process $\{\xi_s,\zeta_s\}$ is
  trivial.
\end{lem}

\textbf{Proof:}
  Define $\mathcal{S}_t = \sigma(\xi_s,s\geq t)$, $\mathcal{T}_t =
  \sigma(\zeta_s, s\geq t)$ and $\mathcal{U}_t = \sigma(\xi_s,
  \zeta_s, s\geq t)$. As pointed out by Ibragimov and Linnik (1971, p. 303) (for
  regularity), to prove that $\mathcal{U}_\infty$ is trivial, it
  suffices to prove that for all $\mathcal{U}_0$-measurable zero
  mean
  random variables $\eta$ such that $\mathbf E [\eta^2] \leq 1$,
  $\mathbf E[\eta \mid \mathcal{U}_t]$ converges to 0 in quadratic
  mean. By standard arguments, it suffices to prove this for a random
  variable $\eta$ that can be expressed as $\eta = \eta_1\eta_2$ with
  $\eta_1$ $\mathcal S_0$-measurable and $\eta_2$ $\mathcal
  T_0$-measurable and, without loss of generality, both with zero mean.
  Then, by independence of $\{\xi_s\}$ and~$\{\zeta_s\}$,
  $$
    \mathbf E[\eta \mid \mathcal{U}_t] = \mathbf E[\eta_1 \mid
    \mathcal{S}_t] \times \mathbf E[\eta_2 \mid \mathcal{T}_t] \; .
  $$
  Since $\mathcal S_\infty$ and $\mathcal T_\infty$ are trivial, both
  terms in the right hand side above tend to 0 in q.m. By
  independence, their product also tends to 0 in q.m.
$\Box$

\end{document}